\documentclass[11pt]{article}
\usepackage{amssymb,latexsym,color,amsmath,pifont,epsfig,cite}
\usepackage[top=0.8in,bottom=1in,left=1in,right=1in]{geometry}
 \usepackage{graphicx}   
\usepackage{epsfig}
\usepackage{graphicx,graphics}
\usepackage{cite}
\usepackage{mathrsfs}
\usepackage{times}
\usepackage{mathptmx}
\usepackage{amssymb,amsmath,amsthm}
\usepackage{amsfonts}
\usepackage{txfonts}
\usepackage{subfigure}
\usepackage{amsfonts} \usepackage{txfonts}
\newtheorem{theorem}{\textbf {Theorem}}

\newtheorem{remark}{\textbf{Remark}}
\newtheorem{example}{\textbf{Example}}
\newtheorem{assumption}{\textbf{Assumption}}

\begin{document}
\title{An Adaptive Algorithm for Synchronization in Diffusively Coupled Systems}
\author{S. Yusef Shafi and Murat Arcak\\ Department of Electrical Engineering and Computer Sciences \\University of California, Berkeley\\ \{yusef,arcak\}@eecs.berkeley.edu}
\date{\today}
\maketitle

\begin{abstract}
We present an adaptive algorithm that guarantees synchronization in diffusively coupled systems. We first consider compartmental systems of ODEs, where each compartment represents a spatial domain of components interconnected through diffusion terms with like components in different compartments. Each set of like components may have its own weighted undirected graph describing the topology of the interconnection between compartments. The link weights are updated adaptively according to the magnitude of the difference between neighboring agents connected by the link. We next consider reaction-diffusion PDEs with Neumann boundary conditions, and derive an analogous algorithm guaranteeing spatial homogenization of solutions.  We provide a numerical example demonstrating the results.
\end{abstract}

\section{Introduction}

Spatially distributed models with diffusive coupling are crucial to understanding the dynamical behavior of a range of engineering and biological systems. This form of coupling encompasses, among others, feedback laws for coordination of multi-agent systems, electromechanical coupling of synchronous machines in power systems, and local update laws in distributed agreement algorithms. Synchronization of diffusively coupled models is an active and rich research area \cite{Hale}, with applications to multi-agent systems, power systems, oscillator circuits, physiological processes, etc. 

The majority of the literature assumes a static interconnection between the nodes in full state models \cite{arcak11aut,Nijmeijer,weislotine2005,stan2007,russo2009,ScardoviEtAl,pecora1998} or phase variables in phase coupled oscillator models \cite{Kuramoto1,Strogatz,chopraspong2009,Dorfler}. However, recently, the situation where interconnection strengths are adapted according to local synchronization errors has started to attract interest. In \cite{assenza2011}, the authors proposed a phase-coupled oscillator model in which local interactions were reinforced between agents with similar behavior and weakened between agents with divergent behavior, leading to enhanced local synchronization. In \cite{delellis2009}, the authors analyzed synchronization of oscillators and presented an adaptive approach to establish synchrony. In \cite{demetriou2013synchronization}, the author considered synchronization and consensus in linear parabolic distributed systems, and in \cite{demetriou2013adaptation} presented an adaptive algorithm to guarantee state regulation and improve convergence of coupled agents to common transient trajectories. 

In this note, we present an adaptive algorithm to guarantee synchrony in diffusively coupled systems. In our earlier work for static interconnections, we gave a numerically-verifiable condition on the Jacobian of a vector field to guarantee spatial homogeneity in reaction-diffusion PDEs and coupled compartmental systems of ODEs \cite{arcak11aut}, and generalized the condition to heterogeneous diffusion in \cite{shafi2013acc}. Using these results as a starting point, here we first consider compartmental models and derive adaptive laws that update interconnection strengths locally to achieve sufficient connectivity for synchronization. We next consider reaction-diffusion partial differential equations, and show that a similar control law that adapts the strength of diffusion coefficients guarantees spatial homogeneity. We present a numerical example that demonstrates the effectiveness of adaptation in enhancing synchrony and lends insight to understanding the structure and most crucial links of the network. 

Our results make several key contributions differing from the existing literature. In \cite{delellis2009}, the authors presented an adaptive law to establish synchrony across agents in a coupled compartmental system of ODEs. They assumed full-state coupling over a single graph via a vector-valued output function. In contrast, we do not assume full-state coupling, and we further allow multiple input-output channels interconnected according to different graphs. In this case, the link weights for each graph are adjusted with a separate update rule. In addition, we present a PDE analogue of the proposed adaptation. In \cite{demetriou2013adaptation}, the author considered a collection of identical linear spatially distributed systems (e.g., linear parabolic PDEs) coupled by a graph with the goals of state regulation to zero and synchrony across agents. However, nonlinear models and nonequilibrium dynamics are not considered. We study nonlinear models and do not make any assumptions on the attractors of this model. This allows us to achieve synchronization for limit cycle oscillators, multi-stable systems, etc. Furthermore, to our knowledge the literature does not address the question of spatial homogenization in reaction-diffusion PDEs in which the coefficients of the elliptic operator vary in time.

\section{Compartmental ODEs}
Let $\mathcal{G}$ be an undirected, connected graph with $N$ nodes and $M$ links, where the nodes $i=1,\cdots,N$ represent the dynamical systems:
\begin{eqnarray}\label{initial}
\dot{x}_i&=&f(x_i)+B\sum_{j=1}^Nk_{ij}\,(y_j-y_i) \quad i=1,\cdots,N
\\
y_i&=&Cx_i \label{initial2}
\end{eqnarray}
in which $x_i\in \mathbb{R}^n$, $B\in \mathbb{R}^{n\times p}$, and $C\in \mathbb{R}^{p\times n}$, $f(\cdot)$ is a continuously differentiable vector field, and the scalars
$k_{ij}=k_{j\, i}$ for each pair $(i,j)$.  Nodes $i$ and $j$ are called neighbors in $\mathcal{G}$ if there is a link in $\mathcal{G}$ connecting $i$ with $j$. We take $k_{ij}=0$ when nodes $i$ and $j$ are not neighbors in $\mathcal{G}$ so that the dynamical systems defined by  (\ref{initial})-(\ref{initial2}), $i=1,\cdots,N$, are coupled according to the graph structure.
When $i$ and $j$ are neighbors, $k_{ij}=k_{j\, i}$ is updated according to:
\begin{equation}\label{update}
\dot{k}_{ij}=\gamma_{ij}\, (y_i-y_j)^T(y_i-y_j)
\end{equation}
where $\gamma_{ij}=\gamma_{j\,i}>0$ is an adaptation gain to be selected by the designer. Thus, there are $M$  independent variables updated as in (\ref{update}), corresponding to each link of the graph.

Define
\begin{equation}\label{defbar}
\bar{x}:=\frac{1}{N}(x_1+\cdots+x_N),
\quad
\tilde{x}_i:=x_i-\bar{x},
\quad \mbox{and} \quad
\tilde{y}_i:=C\tilde{x}_i.
\end{equation}
To guarantee that
$\tilde{y}_i(t) \rightarrow 0$ as $t\rightarrow \infty$, that is, the outputs of the dynamical systems synchronize, we restrict the matrices $B$, $C$, and the Jacobian
$$
J(x)=\frac{\partial f(x)}{\partial x}
$$
with the following assumption:
\begin{assumption}\label{OFP}
There exist a convex set $\mathcal{X}\subset \mathbb{R}^n$, a constant $\theta>0$, and a matrix $P=P^T>0$ such that:
\begin{eqnarray}\label{OFP1}
&& PJ(x)+J(x)^TP\le \theta C^TC \quad \forall x\in \mathcal{X}\\
&& PB=C^T\label{OFP2}.
\end{eqnarray}
\end{assumption}

\begin{theorem}\label{odethm}
Consider the interconnected system (\ref{initial})-(\ref{initial2}), $i=1,\cdots,N$, where $k_{ij}=k_{j\, i}$ is  updated according to (\ref{update}) when nodes $i$ and $j$ are neighbors in $\mathcal{G}$ and is interpreted as zero otherwise, and suppose Assumption \ref{OFP} holds.
If the solutions are bounded and $x_i(t)\in \mathcal{X}$ for all $t\ge 0$, $i=1,\cdots,N$, then $\tilde{y}_i(t) \rightarrow 0$ as $t\rightarrow \infty$.
\hfill $\Box$
\end{theorem}

\noindent
{\bf Proof of Theorem \ref{odethm}:}
We define $\tilde{k}_{ij}={k}_{ij}-{k}_{ij}^*$ where
\begin{equation}\label{kdef}
k^*_{ij}=\left\{ \begin{array}{ll} k^* & \mbox{if $i$ and $j$ are neighbors in $\mathcal{G}$}\\
0 & \mbox{otherwise} \end{array}\right.
\end{equation}
and $k^*$ is a constant to be selected. We then introduce the Lyapunov function:
\begin{equation}
V=\sum_{i=1}^N \tilde{x}_i^TP\tilde{x}_i+\sum_{i=1}^N \sum_{j=1}^N\frac{1}{2\gamma_{ij}}\tilde{k}_{ij}^2
\end{equation}
where (\ref{kdef}) implies that $\tilde{k}_{ij}={k}_{ij}=0$ for pairs $(i,j)$ that are not neighbors in $\mathcal{G}$.
Taking the derivative of $V$ with respect to time and substituting (\ref{initial}) and (\ref{update}), we get:
\begin{equation}
\dot{V}=\sum_{i=1}^N 2\tilde{x}_i^TP\left(f(x_i)-\dot{\bar{x}}+B\sum_{j=1}^Nk_{ij}\,(y_j-y_i)\right)+\sum_{i=1}^N \sum_{j=1}^N\tilde{k}_{ij}(y_i-y_j)^T(y_i-y_j).
\end{equation}
We then substitute $y_i-y_j=\tilde{y}_i-\tilde{y}_j$ and $\tilde{x}_i^TPB=\tilde{y}_i^T$, which follows from (\ref{OFP2}), and obtain:
\begin{equation}\label{dot2}
\dot{V}=\sum_{i=1}^N 2\tilde{x}_i^TP(f(x_i)-\dot{\bar{x}}))+2\sum_{i=1}^N \sum_{j=1}^Nk_{ij}\,\tilde{y}_i^T(\tilde{y}_j-\tilde{y}_i)+\sum_{i=1}^N \sum_{j=1}^N\tilde{k}_{ij}(\tilde{y}_i-\tilde{y}_j)^T(\tilde{y}_i-\tilde{y}_j).
\end{equation}
Next, we note from (\ref{defbar}) that $\sum_{i=1}^N \tilde{x}_i=0$, and add
\begin{equation}
\sum_{i=1}^N 2\tilde{x}_i^TP(\dot{\bar{x}}-f(\bar{x})))=2\left(\sum_{i=1}^N \tilde{x}_i\right)^TP(\dot{\bar{x}}-f(\bar{x})))=0
\end{equation}
to the right-hand side of (\ref{dot2}):
\begin{equation}\label{dot3}
\dot{V}=\sum_{i=1}^N 2\tilde{x}_i^TP(f(x_i)-f(\bar{x})))+2\sum_{i=1}^N \sum_{j=1}^Nk_{ij}\,\tilde{y}_i^T(\tilde{y}_j-\tilde{y}_i)+\sum_{i=1}^N \sum_{j=1}^N\tilde{k}_{ij}\tilde{k}_{ij}(\tilde{y}_i-\tilde{y}_j)^T(\tilde{y}_i-\tilde{y}_j).
\end{equation}
Since
\begin{equation}
f(x_i)-f(\bar{x})=\int_0^1 J(\bar{x}+s\tilde{x}_i)\tilde{x}_i\,ds
\end{equation}
by the Mean Value Theorem, inequality (\ref{OFP1}) yields:
\begin{equation}\label{this}
\tilde{x}_i^TP(f(x_i)-f(\bar{x}))
=\frac{1}{2} \tilde{x}_i^T\left(\int_0^1(PJ(\bar{x}+s\tilde{x}_i)+J^T(\bar{x}+s\tilde{x}_i)P)ds\right)\tilde{x}_i
\le \frac{\theta}{2}\tilde{x}_i^TC^TC\tilde{x}_i=\frac{\theta}{2}\tilde{y}_i^T\tilde{y}_i,
\end{equation}
and substitution of (\ref{this}) in (\ref{dot3}) gives:
\begin{equation}\label{dot4}
\dot{V}\le \theta \sum_{i=1}^{N} \tilde{y}_i^T\tilde{y}_i+2\sum_{i=1}^N \sum_{j=1}^Nk_{ij}\,\tilde{y}_i^T(\tilde{y}_j-\tilde{y}_i)+\sum_{i=1}^N \sum_{j=1}^N\tilde{k}_{ij}(\tilde{y}_i-\tilde{y}_j)^T(\tilde{y}_i-\tilde{y}_j).
\end{equation}
To further simplify (\ref{dot4}), we note that
\begin{equation}
\sum_{i=1}^N \sum_{j=1}^Nk_{ij}\,\tilde{y}_i^T(\tilde{y}_j-\tilde{y}_i)=\sum_{i=1}^N \sum_{j=1}^Nk_{ij}\,\tilde{y}_j^T(\tilde{y}_i-\tilde{y}_j)
\end{equation}
which follows by swapping the indices $i$ and $j$ and substituting $k_{ij}=k_{j\, i}$. Thus,
\begin{equation}
2\sum_{i=1}^N \sum_{j=1}^Nk_{ij}\,\tilde{y}_i^T(\tilde{y}_j-\tilde{y}_i)=\sum_{i=1}^N \sum_{j=1}^Nk_{ij}\,\tilde{y}_i^T(\tilde{y}_j-\tilde{y}_i)-\sum_{i=1}^N \sum_{j=1}^Nk_{ij}\,\tilde{y}_j^T(\tilde{y}_j-\tilde{y}_i)=-\sum_{i=1}^N \sum_{j=1}^N{k}_{ij}(\tilde{y}_i-\tilde{y}_j)^T(\tilde{y}_i-\tilde{y}_j),
\end{equation}
and (\ref{dot4}) becomes:
\begin{equation}\label{dot5}
\dot{V}\le \theta \sum_{i=1}^{N} \tilde{y}_i^T\tilde{y}_i-\sum_{i=1}^N \sum_{j=1}^N{k}^*_{ij}(\tilde{y}_i-\tilde{y}_j)^T(\tilde{y}_i-\tilde{y}_j).
\end{equation}

Next, we assign an arbitrary orientation to the links of the graph $\mathcal{G}$, label the links $\ell=1,\cdots,M$, and introduce the $N\times M$ incidence matrix:
\begin{equation}
E_{i\ell}=\left\{ \begin{array}{ll}
1 & \mbox{if node $i$ is the head of link $\ell$}\\
-1 & \mbox{if node $i$ is the tail of link $\ell$}\\
0 & \mbox{if node $i$ is not connected to link $\ell$.}\\
\end{array}\right.
\end{equation}
Defining $\tilde{Y}=[\tilde{y}_1\ \cdots \tilde{y}_N]^T$, we note that $(E\otimes I_p)^T\tilde{Y}$ is a column vector which is a concatenation of $p$-dimensional components and the $\ell$th such component is $\tilde{y}_i-\tilde{y}_j$ where $i$ is the head and $j$ is the tail of link $\ell$.
It then follows from (\ref{kdef}) that:
\begin{equation}
\sum_{i=1}^N \sum_{j=1}^N{k}^*_{ij}(\tilde{y}_i-\tilde{y}_j)^T(\tilde{y}_i-\tilde{y}_j)=2k^*\tilde{Y}^T(E\otimes I_p)(E\otimes I_p)^T\tilde{Y}=2k^*\tilde{Y}^T(EE^T\otimes I_p)\tilde{Y}.
\end{equation}
Since $EE^T$ is the Laplacian matrix for the graph $\mathcal{G}$, its smallest eigenvalue is $\lambda_1=0$ and the vector of ones $\mathbf{1}_N$ is a corresponding eigenvector. Likewise, for $EE^T\otimes I_p$, $\lambda_1=0$ has multiplicity $p$ and the corresponding eigenspace is the range of $\mathbf{1}_N\otimes I_p$.
Because $\mathcal{G}$ is connected, the second smallest eigenvalue $\lambda_2$ is strictly positive and, since $\tilde{Y}^T(\mathbf{1}_N\otimes I_p)=0$ from (\ref{defbar}), the following inequality holds:
\begin{equation}
\tilde{Y}^T(EE^T\otimes I_p)\tilde{Y}\ge \lambda_2 \tilde{Y}^T\tilde{Y}.
\end{equation}
Thus,
\begin{equation}\label{dot6}
\dot{V}\le -(2k^*\lambda_2-\theta)\tilde{Y}^T\tilde{Y}
\end{equation}
and choosing $k^*$ large enough that $\epsilon:=2k^*\lambda_2-\theta>0$ guarantees:
\begin{equation}\label{dot7}
\dot{V}\le -\epsilon \tilde{Y}^T\tilde{Y}.
\end{equation}
By integrating both sides of the inequality (\ref{dot7}), we conclude that  $\tilde{y}_i(t)$ is in $\mathcal{L}_2$,  $i=1,\cdots,N$. Furthermore, the boundedness of solutions implies that $\dot{x}_i(t)$ and, thus $\dot{\tilde{y}}_i(t)$ is bounded. Barbalat's Lemma \cite{khalil} then guarantees $\tilde{y}_i(t) \rightarrow 0$ as $t\rightarrow \infty$. \hfill $\Box$
\smallskip

\begin{remark}
An extension of Theorem \ref{odethm} to the case of multiple input-output channels, connected according to different graphs, is straightforward.  The system now takes the form: 
\begin{eqnarray}\label{initial3}
\dot{x}_i&=&f(x_i)+\sum_{q=1}^mB^{(q)}\sum_{j=1}^Nk_{ij}^{(q)}\,(y^{(q)}_j-y^{(q)}_i) \quad i=1,\cdots,N
\\
y_i^{(q)}&=&C^{(q)}x_i \label{initial4}
\end{eqnarray}
where $B^{(q)}\in \mathbb{R}^{n\times p_q}$ and $C^{(q)}\in \mathbb{R}^{p_q\times n}$, $q=1,\cdots,m$. A graph $\mathcal{G}^{(q)}$ is defined for each channel $q$ and $k_{ij}^{(q)}=k_{j\, i}^{(q)}\neq 0$ only when nodes $i$ and $j$ are adjacent in $\mathcal{G}^{(q)}$.
 The update rule (\ref{update}) then becomes:
\begin{equation}\label{MIMOupdate}
\dot{k}^{(q)}_{ij}=\gamma^{(q)}_{ij}\, (y^{(q)}_i-y^{(q)}_j)^T(y^{(q)}_i-y^{(q)}_j), \quad  \gamma^{(q)}_{ij}>0,  \quad q=1,\cdots,m.
\end{equation}
To prove synchronization we now ask that (\ref{OFP1})-(\ref{OFP2}) in Assumption \ref{OFP} hold for the matrices $B=[B^{(1)} \cdots B^{(m)}]$ and $C^T=[{C^{(1)}}^T \cdots {C^{(m)}}^T]$. 
In fact, one can relax (\ref{OFP2}) as:
\begin{equation}\label{OFP2b}
PB=[\omega^{(1)}{C^{(1)}}^T \cdots \omega^{(m)}{C^{(m)}}^T],
\end{equation}
where $\omega^{(q)}>0$, $q=1,\cdots,m$.  To accommodate the ``multipliers" $\omega^{(q)}$, the Lyapunov function in the proof of Theorem \ref{odethm} is modified as:
\begin{equation}
V=\sum_{i=1}^N \tilde{x}_i^TP\tilde{x}_i+\sum_{q=1}^m\sum_{i=1}^N \sum_{j=1}^N\frac{\omega^{(q)}}{2\gamma^{(q)}_{ij}}(\tilde{k}_{ij}^{(q)})^2.
\end{equation}
The steps of the proof are otherwise identical and are not repeated to avoid excessive notation.
\end{remark}

\begin{remark}\label{boundedness}
Since the proof above analyzes the evolution of $x_i$ relative to the average $\bar{x}$,  it cannot reach any conclusions about the absolute behavior of the variables $x_i$.
Thus, boundedness of the solutions does not follow from the proof and was assumed in the theorem. However, it is possible to conclude boundedness with an additional restriction on the vector field $f(x)$:
Since $k_{ij}=k_{ji}$, the coupling terms in (\ref{initial}) do not affect the evolution of the average $\bar{x}$, which is governed by:
\begin{equation}\label{BIBS}
\dot{\bar{x}}=\frac{1}{N}\sum_{i=1}^Nf(x_i)=\frac{1}{N}\sum_{i=1}^Nf(\bar{x}+\tilde{x}_i).
\end{equation}
If this system has a bounded-input-bounded-state (BIBS) property when $\tilde{x}_i$ are interpreted as inputs, then we conclude boundedness of all solutions.
This follows because the Lyapunov arguments in the proof show that $\tilde{x}_i(t)$ are bounded on the maximal interval of existence $[0,t_f)$ with bounds that do not depend on $t_f$ and, thus, a similar conclusion holds for $\bar{x}(t)$ implying that $t_f=\infty$. The Lyapunov function then establishes boundedness of $\tilde{x}_i(t)$ and $k_{ij}(t)$, and the BIBS property above guarantees that all solutions are bounded.

A further assumption of the theorem is that $x_i(t)\in \mathcal{X}$ for all $t\ge 0$.  Thus, when the set $\mathcal{X}$ where (\ref{OFP1}) holds is a strict subset of $\mathbb{R}^n$, we have to independently show that $x_i(t)$ remains in $\mathcal{X}$. One can do this by establishing the invariance of the set $\mathcal{X}^N\times \mathbb{R}^M$ for (\ref{initial})-(\ref{update}). If $\mathcal{X}^N\times \mathbb{R}^M$ is not invariant, then an appropriate reachability analysis can be used to identify a set of initial conditions such that the trajectories starting in this set do not leave $\mathcal{X}^N\times \mathbb{R}^M$. \hfill $\Box$
\end{remark}

\begin{figure}[t]
\vspace{-3.0cm}
\begin{center}
\mbox{}\setlength{\unitlength}{.9mm}
\begin{picture}(60,70)
\put(-10,0){\psfig{figure=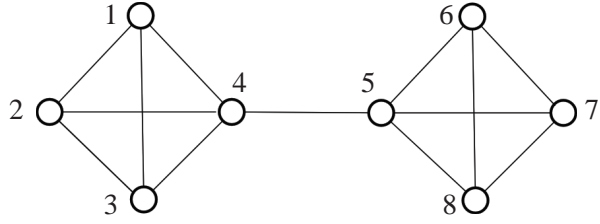,width=80\unitlength}}
\put(-14,14){$2$}\put(19,18){$4$}\put(38,18){$5$}\put(71,14){$7$}
\put(0,0){$3$}\put(50,0){$8$}
\put(0,28){$1$}\put(49.50,28){$6$}
\end{picture}
\vspace{-.2cm}
 \caption{\small An eight-node ``barbell" graph.}  \label{barbell}
\end{center}
\vspace{-.6cm}
\end{figure}

\begin{figure}[h!]
\vspace{1.9cm}
\begin{center}
\mbox{}\setlength{\unitlength}{.7mm}
\begin{picture}(60,70)
\put(-53,35){\psfig{figure=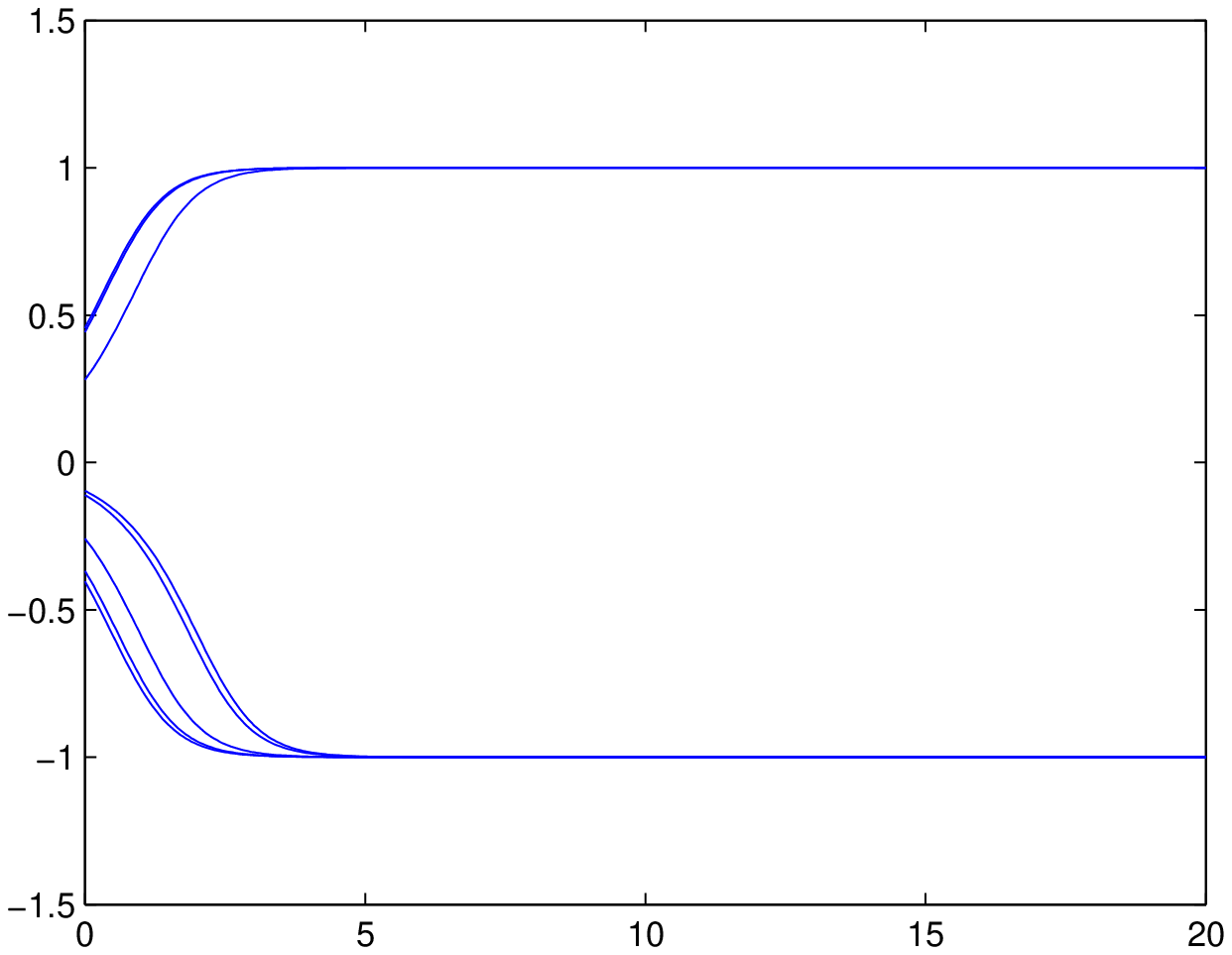,width=80\unitlength}}
\put(-57,87){(A)}\put(-58,65){$x_{i}(t)$}\put(22,41){$t$}
\put(32,35){\psfig{figure=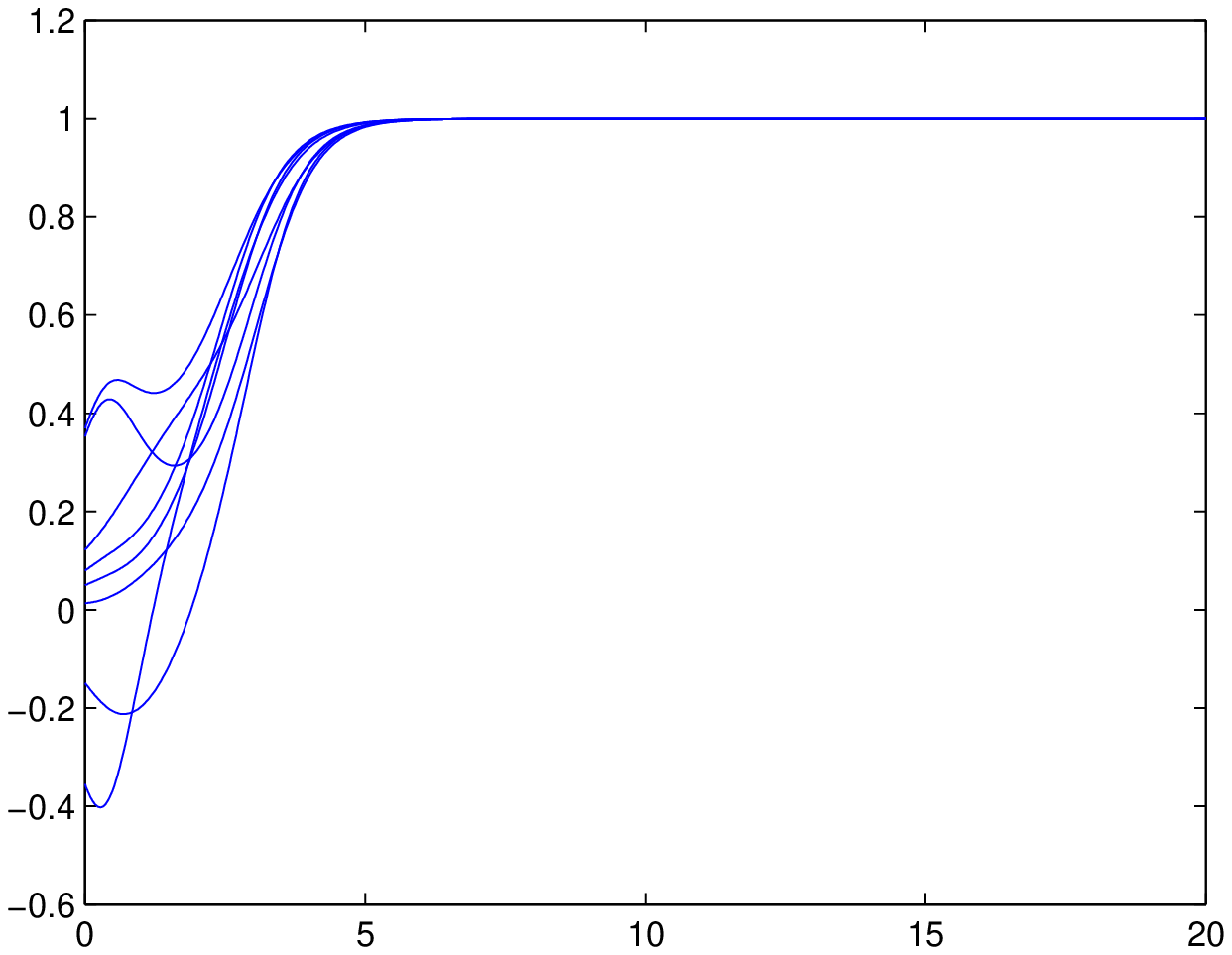,width=80\unitlength}}
\put(28,87){(B)}\put(27,65){$x_{i}(t)$}\put(107,41){$t$}
\put(-53,-23){\psfig{figure=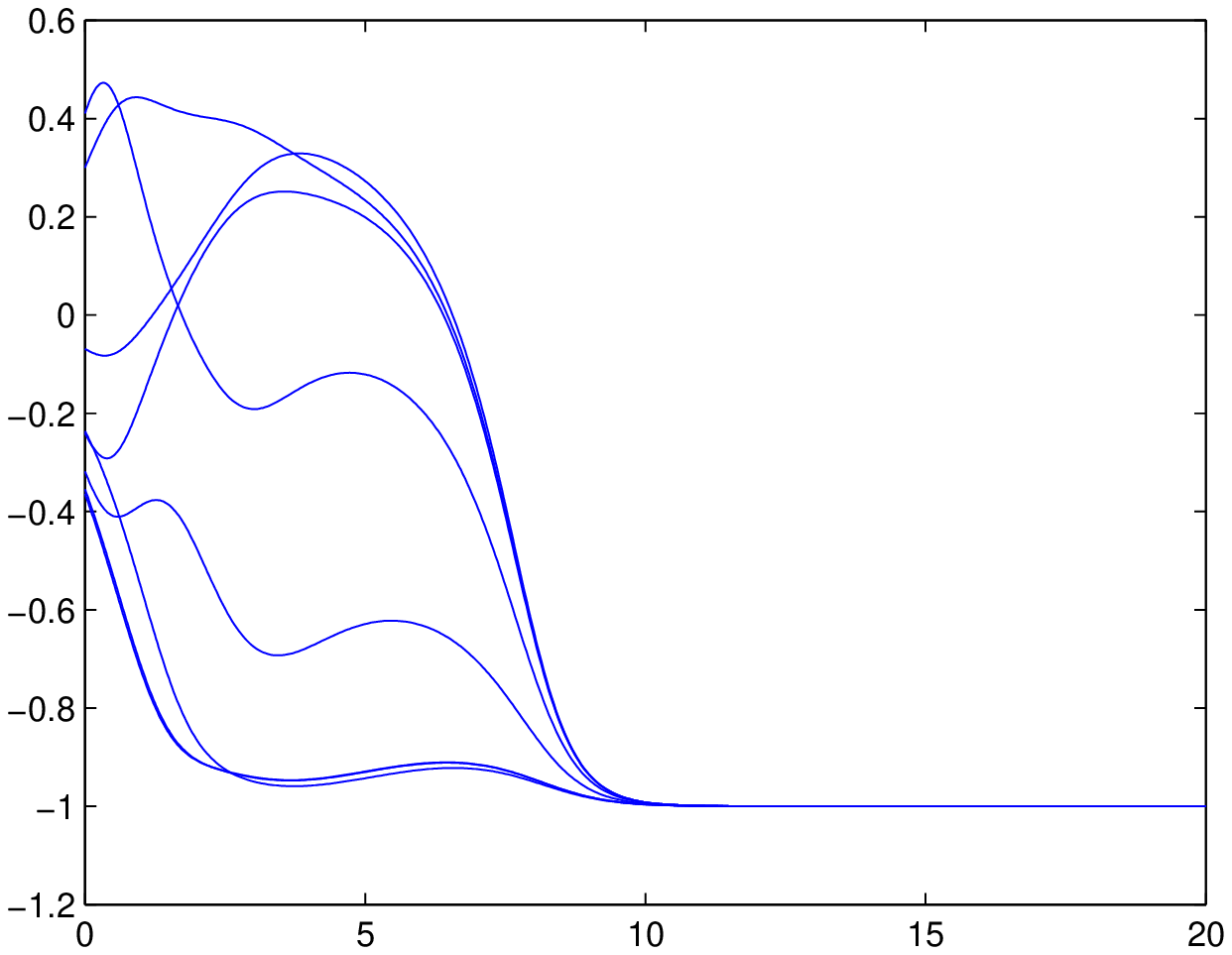,width=80\unitlength}}
\put(-57,28){(C)}\put(-58,6){$x_{i}(t)$}\put(22,-17){$t$}
\put(32,-23){\psfig{figure=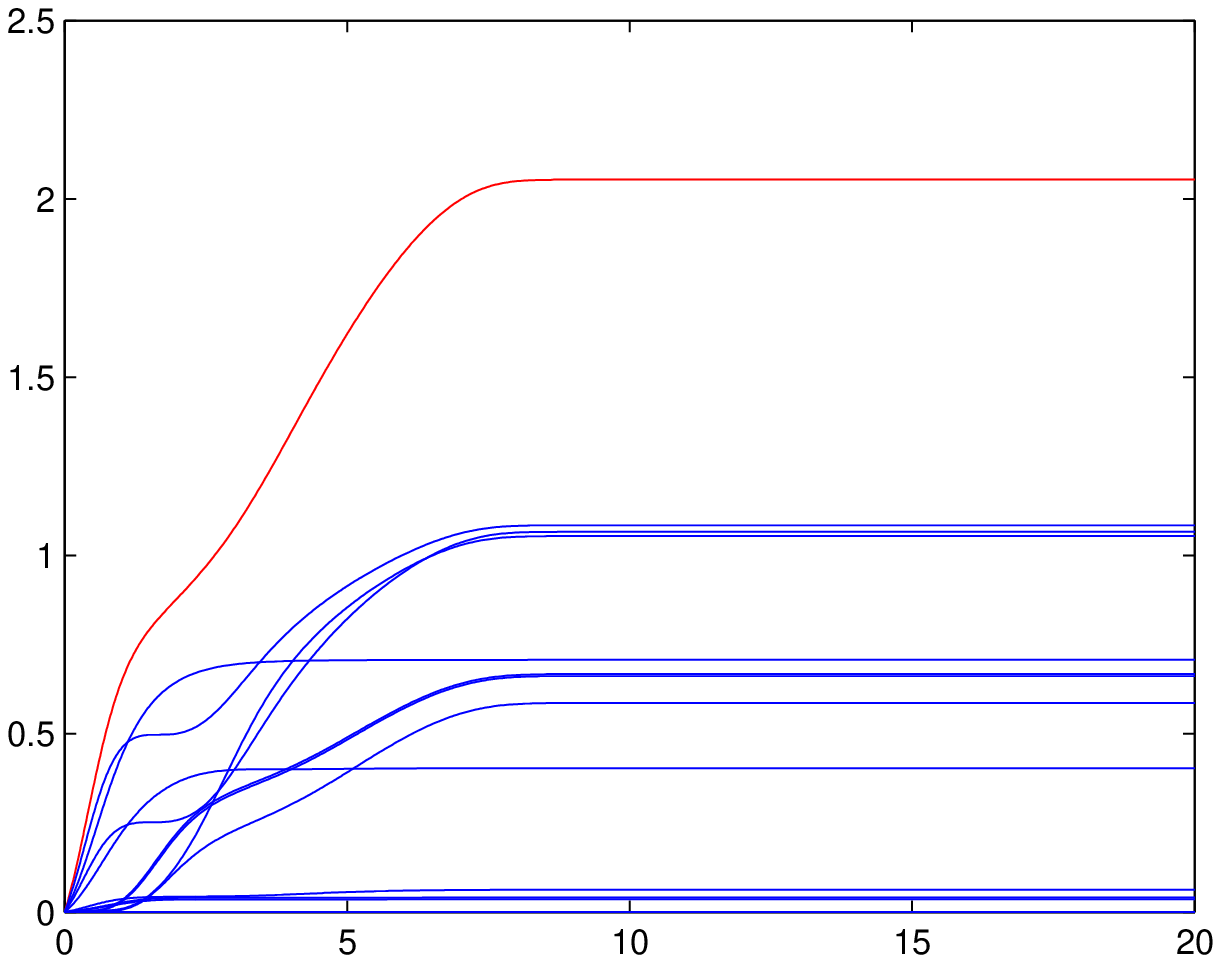,width=80\unitlength}}
\put(28,28){(D)}\put(28,6){$k_{ij}(t)$}\put(107,-17){$t$}
\end{picture}
\vspace{1.1cm}
 \caption{\small Simulations for the graph in Figure \ref{barbell} where each node is governed by (\ref{initial})-(\ref{initial2}) with $B=C=1$ and
$f(x)=x-x^3$. (A) When the link weights are set to zero and the adaptation is turned off, $x_i(t)$ converge to one of the two stable equilibria  $\pm 1$ and do not synchronize. (B)-(C) When the adaptation is turned on, $x_i(t)$ synchronize and converge to identical states. (D) The evolution of link weights $k_{ij}(t)$ in the simulation corresponding to Figure \ref{ex}C. The red curve is $k_{45}(t)$ for the bottleneck link $(4,5)$.}  \label{ex}
\end{center}
\vspace{-.6cm}
\end{figure}

\begin{example} Consider the graph in Figure \ref{barbell} and supposed the nodes are governed by (\ref{initial})-(\ref{initial2}) where $B=C=1$ and
\begin{equation}\label{bistable}f(x)=x-x^3.
  \end{equation}
Thus, Assumption \ref{OFP} holds in $\mathcal{X}=\mathbb{R}$ with $P=1$ and $\theta=2$. The boundedness of solutions condition in Theorem \ref{odethm} follows from the BIBS property of (\ref{BIBS}), as in Remark \ref{boundedness}. To see this BIBS property, note from (\ref{bistable}) that $|\bar{x}|>1+\max_i\{|\tilde{x}_i|\}$ implies $\bar{x}f(\bar{x}+\tilde{x}_i)<0$. This guarantees that the solutions of (\ref{BIBS}) satisfy $|\bar{x}(t)|\le \max\{|\bar{x}(0)|,1+\max_i \sup_{t\ge 0}|\tilde{x}_i(t)|\}$ for all $t\ge 0$.

Note from (\ref{bistable}) that each node is a bistable system with stable equilibria at $x_i=\pm 1$ and a saddle point at $x_i=0$.  Indeed, when we set the link weights to zero and turn off the adaptation, each $x_i(t)$ evolves independently and converges to $+1$ or $-1$ (Figure \ref{ex}A). Next, we turn on the adaptation with
gain $\gamma_{ij}=1$ and  initial condition $k_{ij}(0)=0$ for each link. Figures \ref{ex}B and \ref{ex}C confirm that the nodes now synchronize, converging to a consensus value of $+1$ or $-1$.
Note from (\ref{update}) that the final value reached by each $k_{ij}(t)$  is the squared $L_2$ norm for the synchronization error $x_i(t)-x_j(t)$. As one may expect from the structure of the ``barbell" graph in Figure \ref{barbell}, the red curve in Figure \ref{ex}D corresponding to the bottleneck link $(4,5)$ reached the largest value, indicating a high ``stress" on this link.
\hfill \end{example}

\section{Reaction-Diffusion PDEs}\label{pdesec}
Let $\Omega$ be a bounded and connected domain in $\mathbb{R}^r$ with smooth boundary $\partial \Omega$, and
consider the PDE:
\begin{eqnarray}\label{rdnet}
\frac{\partial x(t,\xi)}{\partial t}&=&f(x(t,\xi))+\sum_{\ell=1}^{p}B_{\ell} \nabla\cdot (k(t,\xi)\nabla y_{\ell}(t,\xi)),\\
\label{rdout}
y_{\ell}(t,\xi)&=&C_{\ell} x(t,\xi)
\end{eqnarray}
where $\xi \in \Omega$ is the spatial variable, $x(t,\xi)\in \mathbb{R}^n$, $k(t,\xi)\in \mathbb{R}$, $f(\cdot)$ is a continuously differentiable vector field, $B=[B_1\cdots B_p]\in \mathbb{R}^{n\times p}$, $C^T=[C_1^T \cdots C_p^T]\in \mathbb{R}^{n\times p}$, $\nabla \cdot$ is the divergence operator and $\nabla$ represents the gradient with respect to the spatial variable $\xi$.  
We assume Neumann boundary conditions:
\begin{equation}
\nabla x_i(t,\xi) \cdot \hat{n}(\xi)=0 \quad \forall \xi \in \partial \Omega, \ \forall t\ge 0, \quad i=1,\cdots,n \label{bc}
\end{equation}
where ``$\cdot$" is the inner product in $\mathbb{R}^r$,  $x_i(t,\xi)$ denotes the $i$th entry of the vector $x(t,\xi)$ and $\hat{n}(\xi)$ is a vector normal to the boundary $\partial \Omega$.

In analogy with (\ref{update}), we introduce the update law:
\begin{equation}\label{rdupdate}
\frac{\partial k(t,\xi)}{\partial t}=\gamma(\xi)\sum_{\ell=1}^{p}\nabla y_\ell(t,\xi)\cdot \nabla y_\ell(t,\xi),
\end{equation}
where $\gamma(\xi)>0$ is a design choice.
Define:
\begin{equation}\label{pidef}
\bar{x}(t):=\frac{1}{|\Omega|}\int_\Omega {x}(t,\xi) d\xi,
\quad
\tilde{x}(t,\xi):={x}(t,\xi)-\bar{x}(t),
\quad
\tilde{y}(t,\xi):=C\tilde{x}(t,\xi).
\end{equation}
In Theorem \ref{pdethm} below, we give conditions that guarantee the following output synchronization property:
\begin{equation}\label{rdsync}
\lim_{t\rightarrow \infty}\int_\Omega |\tilde{y}(t,\xi)|^2d\xi =0
\end{equation}
where $|\cdot |$ denotes the Euclidean norm.
\begin{theorem}\label{pdethm}
Consider the system (\ref{rdnet})-(\ref{rdout}) with boundary condition (\ref{bc}), and suppose Assumption \ref{OFP} holds. Then, the update law (\ref{rdupdate}) guarantees (\ref{rdsync}) for every bounded classical solution that satisfies
$x(t,\xi)\in \mathcal{X}$ for all $t\ge 0$.
\hfill $\Box$
\end{theorem}

\noindent
{\bf Proof of Theorem \ref{pdethm}:}
Define:
\begin{equation}
V(t)=\int_\Omega \tilde{x}^T(t,\xi)P\tilde{x}(t,\xi)d\xi+\int_\Omega \frac{1}{\gamma(\xi)}|\tilde{k}(t,\xi)|^2d\xi
\end{equation}
where $\tilde{k}(t,\xi)={k}(t,\xi)-k^*$, and $k^*$ is to be selected. Taking derivatives with respect to time and using (\ref{rdupdate}), we get:
\begin{equation}\label{rdot1}
\dot{V}(t)=2\int_\Omega \tilde{x}^T(t,\xi)P\frac{\partial \tilde{x}(t,\xi)}{\partial t}d\xi+2\sum_{\ell=1}^p\int_\Omega \tilde{k}(t,\xi)\nabla y_\ell(t,\xi)\cdot \nabla y_\ell(t,\xi)d\xi.
\end{equation}
It follows from (\ref{rdnet}) that:
\begin{equation}\label{rdot2}
\dot{V}(t)=2\int_\Omega \tilde{x}^T(t,\xi)P\left(f(x(t,\xi))-\dot{\bar{x}}+\sum_{\ell=1}^pB_\ell\nabla\cdot (k(t,\xi)\nabla y_\ell(t,\xi))\right)d\xi+2\sum_{\ell=1}^p\int_\Omega \tilde{k}(t,\xi)\nabla y_\ell(t,\xi)\cdot \nabla y_\ell(t,\xi)d\xi.
\end{equation}
Next, substituting $\nabla y_\ell(t,\xi)=\nabla \tilde{y}_\ell(t,\xi)$ and $\tilde{x}^T(t,\xi)PB_\ell=\tilde{y}_\ell(t,\xi)$, which follows from (\ref{OFP2}), we obtain:
\begin{eqnarray}\label{rdot3}
\dot{V}(t)&=&2\int_\Omega \tilde{x}^T(t,\xi)P(f(x(t,\xi))-\dot{\bar{x}})d\xi+2\sum_{\ell=1}^p\int_\Omega \tilde{y}_\ell(t,\xi)\nabla\cdot (k(t,\xi)\nabla \tilde{y}_\ell(t,\xi)))d\xi\\&& +2\sum_{\ell=1}^p\int_\Omega \tilde{k}(t,\xi)\nabla \tilde{y}_\ell(t,\xi)\cdot \nabla \tilde{y}_\ell(t,\xi)d\xi.\nonumber
\end{eqnarray}
Since
\begin{equation}
\int_\Omega \tilde{x}^T(t,\xi)P(\dot{\bar{x}}(t)-f(\bar{x}(t)))d\xi=\left(\int_\Omega \tilde{x}^T(t,\xi)d\xi\right)P(\dot{\bar{x}}(t)-f(\bar{x}(t)))=0,
\end{equation}
which follows from (\ref{pidef}), we rewrite (\ref{rdot3}) as:
\begin{eqnarray}\label{rdot4}
\dot{V}(t)&=&2\int_\Omega \tilde{x}^T(t,\xi)P(f(x(t,\xi))-f(\bar{x}(t)))d\xi+2\sum_{\ell=1}^p\int_\Omega \tilde{y}_\ell(t,\xi)\nabla\cdot (k(t,\xi)\nabla \tilde{y}_\ell(t,\xi)))d\xi\\&& +2\sum_{\ell=1}^p\int_\Omega \tilde{k}(t,\xi)\nabla \tilde{y}_\ell(t,\xi)\cdot \nabla \tilde{y}_\ell(t,\xi)d\xi.\nonumber
\end{eqnarray}
It then follows from (\ref{this}) that
\begin{equation}\label{rdot5}
\dot{V}(t)\le {\theta}\int_\Omega |\tilde{y}(t,\xi)|^2d\xi+2\sum_{\ell=1}^p\int_\Omega \tilde{y}_\ell(t,\xi)\nabla\cdot (k(t,\xi)\nabla \tilde{y}_\ell(t,\xi)))d\xi +2\sum_{\ell=1}^p\int_\Omega \tilde{k}(t,\xi)\nabla \tilde{y}_\ell(t,\xi)\cdot \nabla \tilde{y}_\ell(t,\xi)d\xi.
\end{equation}
We now claim that:
\begin{equation}\label{claim}
\int_\Omega \tilde{y}_\ell(t,\xi)\nabla\cdot (k(t,\xi)\nabla \tilde{y}_\ell(t,\xi)))d\xi=-\int_\Omega {k}(t,\xi)\nabla \tilde{y}_\ell(t,\xi)\cdot \nabla \tilde{y}_\ell(t,\xi)d\xi.
\end{equation}
This follows by first applying the identity $\nabla \cdot (fF)=f \nabla \cdot F+F\cdot \nabla f$, which holds when $f$ is scalar valued, with $F=k(t,\xi)\nabla \tilde{y}_\ell(t,\xi)$ and $f=\tilde{y}_\ell(t,\xi)$, next integrating both sides of the identity over $\Omega$, and finally noting that the left-hand side is zero, since:
\begin{equation}
\int_\Omega \nabla \cdot (\tilde{y}_\ell(t,\xi)k(t,\xi)\nabla \tilde{y}_\ell(t,\xi))d\xi=\int_{\partial \Omega} \tilde{y}_\ell(t,\xi)k(t,\xi)\nabla \tilde{y}_\ell(t,\xi) \cdot \hat{n}(\xi)dS
\end{equation}
from the Divergence Theorem and $\nabla \tilde{y}_\ell(t,\xi) \cdot \hat{n}(\xi)=0$ for $\xi \in \partial \Omega$ from the boundary condition (\ref{bc}). Substituting (\ref{claim}) in (\ref{rdot5}), we get:
\begin{equation}\label{rdot6}
\dot{V}(t)\le {\theta}\int_\Omega |\tilde{y}(t,\xi)|^2d\xi-2k^*\sum_{\ell=1}^{p}\int_\Omega \nabla \tilde{y}_\ell(t,\xi)\cdot \nabla \tilde{y}_\ell(t,\xi)d\xi.
\end{equation}
Moreover, because $\int_\Omega \tilde{y}_\ell(t,\xi) d\xi=0$, it follows from the the Poincar\'{e} Inequality \cite[Equation (1.37)]{henrot}
that:
\begin{equation}\label{there}
  \int_\Omega |\nabla{\tilde{y}_\ell(t,\xi)}|^2 d\xi \ge \lambda_2\int_\Omega \tilde{y}_\ell(t,\xi)^2 d\xi
\end{equation}
where $\lambda_2$ denotes the second smallest of the eigenvalues $0=\lambda_1\le \lambda_2 \le \cdots$ of the operator $L=-\nabla^2$ on $\Omega$ with Neumann boundary condition, and $\lambda_2>0$ since $\Omega$ is connected. Thus, (\ref{rdot6}) becomes:
\begin{equation}\label{rdot7}
\dot{V}(t)\le -\left(2k^*\lambda_2-{\theta}\right)\int_\Omega |\tilde{y}(t,\xi)|^2d\xi,
\end{equation}
and choosing $k^*$ large enough that $\epsilon:=2k^*\lambda_2-\theta>0$ guarantees:
\begin{equation}\label{rdot}
\dot{V}(t)\le -\epsilon\int_\Omega |\tilde{y}(t,\xi)|^2d\xi =: -\epsilon W(t).
\end{equation}
This implies that $\lim_{T\rightarrow \infty}\int_0^T W(t)dt$ exists and is bounded. Since $\dot{W}(t)$ is also bounded, it follows from Barbalat's Lemma \cite{khalil} that $W(t)\rightarrow 0$ as $t\rightarrow \infty$ which proves (\ref{rdsync}). \hfill $\Box$

\bibliographystyle{IEEEtran}       
\bibliography{mybib,books,General}  

\begin{thebibliography}{10}
\providecommand{\url}[1]{#1}
\csname url@samestyle\endcsname
\providecommand{\newblock}{\relax}
\providecommand{\bibinfo}[2]{#2}
\providecommand{\BIBentrySTDinterwordspacing}{\spaceskip=0pt\relax}
\providecommand{\BIBentryALTinterwordstretchfactor}{4}
\providecommand{\BIBentryALTinterwordspacing}{\spaceskip=\fontdimen2\font plus
\BIBentryALTinterwordstretchfactor\fontdimen3\font minus
  \fontdimen4\font\relax}
\providecommand{\BIBforeignlanguage}[2]{{%
\expandafter\ifx\csname l@#1\endcsname\relax
\typeout{** WARNING: IEEEtran.bst: No hyphenation pattern has been}%
\typeout{** loaded for the language `#1'. Using the pattern for}%
\typeout{** the default language instead.}%
\else
\language=\csname l@#1\endcsname
\fi
#2}}
\providecommand{\BIBdecl}{\relax}
\BIBdecl

\bibitem{Hale}
J.~Hale, ``Diffusive coupling, dissipation, and synchronization,''
  \emph{Journal of Dynamics and Differential Equations}, vol.~9, no.~1, pp.
  1--52, 1997.

\bibitem{arcak11aut}
M.~Arcak, ``Certifying spatially uniform behavior in reaction-diffusion {PDE}
  and compartmental {ODE} systems,'' \emph{Automatica}, vol.~47, no.~6, pp.
  1219--1229, 2011.

\bibitem{Nijmeijer}
A.~Pogromsky and H.~Nijmeijer, ``Cooperative oscillatory behavior of mutually
  coupled dynamical systems,'' \emph{IEEE Transactions on Circuits and Systems
  I: Fundamental Theory and Applications}, vol.~48, no.~2, pp. 152--162, 2001.

\bibitem{weislotine2005}
W.~Wang and J.-J.~E. Slotine, ``On partial contraction analysis for coupled
  nonlinear oscillators,'' \emph{Biological Cybernetics}, vol.~92, pp. 38--53,
  2005.

\bibitem{stan2007}
G.~Stan and R.~Sepulchre, ``Analysis of interconnected oscillators by
  dissipativity theory,'' \emph{IEEE Transactions on Automatic Control},
  vol.~52, no.~2, pp. 256--270, 2007.

\bibitem{russo2009}
G.~Russo and M.~Di~Bernardo, ``Contraction theory and master stability
  function: linking two approaches to study synchronization of complex
  networks,'' \emph{IEEE Transactions on Circuits and Systems II: Express
  Briefs}, vol.~56, no.~2, pp. 177--181, 2009.

\bibitem{ScardoviEtAl}
L.~Scardovi, M.~Arcak, and E.~Sontag, ``Synchronization of interconnected
  systems with applications to biochemical networks: An input-output
  approach,'' \emph{IEEE Transactions on Automatic Control}, vol.~55, no.~6,
  pp. 1367--1379, 2010.

\bibitem{pecora1998}
L.~M. Pecora and T.~L. Carroll, ``Master stability functions for synchronized
  coupled systems,'' \emph{Physical Review Letters}, vol.~80, no.~10, pp.
  2109--2112, 1998.

\bibitem{Kuramoto1}
Y.~Kuramoto, ``Self-entrainment of a population of coupled non-linear
  oscillators,'' in \emph{International Symposium on Mathematical Problems in
  Theoretical Physics}.\hskip 1em plus 0.5em minus 0.4em\relax Springer, 1975,
  pp. 420--422.

\bibitem{Strogatz}
S.~Strogatz, ``From {K}uramoto to {C}rawford: exploring the onset of
  synchronization in populations of coupled oscillators,'' \emph{Physica D:
  Nonlinear Phenomena}, vol. 143, no.~1, pp. 1--20, 2000.

\bibitem{chopraspong2009}
N.~Chopra and M.~Spong, ``On exponential synchronization of {K}uramoto
  oscillators,'' \emph{IEEE Transactions on Automatic Control}, vol.~54, no.~2,
  pp. 353--357, 2009.

\bibitem{Dorfler}
F.~D{\"o}rfler and F.~Bullo, ``Exploring synchronization in complex oscillator
  networks,'' \emph{arXiv Preprint: arXiv:1209.1335}, 2012.

\bibitem{assenza2011}
S.~Assenza, R.~Guti{\'e}rrez, J.~G{\'o}mez-Garde{\~n}es, V.~Latora, and
  S.~Boccaletti, ``Emergence of structural patterns out of synchronization in
  networks with competitive interactions,'' \emph{Scientific Reports}, vol.~1,
  2011.

\bibitem{delellis2009}
P.~DeLellis, M.~DiBernardo, and F.~Garofalo, ``Novel decentralized adaptive
  strategies for the synchronization of complex networks,'' \emph{Automatica},
  vol.~45, no.~5, pp. 1312--1318, 2009.

\bibitem{demetriou2013synchronization}
M.~A. Demetriou, ``Synchronization and consensus controllers for a class of
  parabolic distributed parameter systems,'' \emph{Systems \& Control Letters},
  vol.~62, no.~1, pp. 70--76, 2013.

\bibitem{demetriou2013adaptation}
------, ``Adaptation and optimization of synchronization gains in networked
  distributed parameter systems,'' \emph{arXiv preprint arXiv:1305.7117}, 2013.

\bibitem{shafi2013acc}
S.~Y. Shafi, Z.~Aminzare, M.~Arcak, and E.~D. Sontag, ``Spatial uniformity in
  diffusively-coupled systems using weighted ${L}_2$ norm contractions,'' in
  \emph{American Control Conference}, 2013, pp. 5619--5624.

\bibitem{khalil}
H.~Khalil, \emph{Nonlinear Systems}.\hskip 1em plus 0.5em minus 0.4em\relax
  Englewood Cliffs, NJ: Prentice Hall, 2002.

\bibitem{henrot}
A.~Henrot, \emph{Extremum problems for eigenvalues of elliptic
  operators}.\hskip 1em plus 0.5em minus 0.4em\relax Basel: Birkhauser, 2006.

\end{thebibliography}

\end{document}